\documentclass[12pt]{amsart}
\usepackage{amsmath,amscd,amsthm,amsfonts, amssymb,amsxtra}
\usepackage[all]{xy}

\newtheorem{thm}{Theorem}[section]  
\newtheorem*{un-no-thm}{Theorem}
\newtheorem{cor}[thm]{Corollary}     
\newtheorem{lem}[thm]{Lemma}         
\newtheorem{prop}[thm]{Proposition}  

\newtheorem{bigthm}{Theorem}
\newtheorem{bigcor}[bigthm]{Corollary}

\theoremstyle{definition}
\newtheorem{defn}[thm]{Definition}   

\theoremstyle{definition}

\theoremstyle{definition}

\theoremstyle{remark}
\newtheorem{rem}[thm]{Remark}
\newtheorem{rems}[thm]{Remarks}

\newtheorem*{acks}{Acknowledgements}
 
\newtheorem*{out}{Outline}

\newtheorem{ex}[thm]{Example}

\begin{document}
\title{Homotopical Intersection Theory, I.}
\date{\today}
\author{John R. Klein}
\address{Wayne State University, Detroit, MI 48202}
\email{klein@math.wayne.edu}
\author{E.\ Bruce Williams}
\address{University of Notre Dame, Notre Dame, IN 46556}
\email{williams.4@nd.edu}
\begin{abstract} We give a new approach to intersection
theory. Our  ``cycles'' are closed manifolds mapping into
compact manifolds and our ``intersections'' are elements
of a homotopy group of a certain Thom space. 
The results are then applied in various contexts,
including fixed point, linking and disjunction problems.
Our main theorems resemble those of Hatcher and Quinn \cite{H-Q},
but our proofs are fundamentally different.
\end{abstract}
\thanks{Both authors are partially supported by the NSF}
\maketitle
\setlength{\parindent}{15pt}
\setlength{\parskip}{1pt plus 0pt minus 1pt}
\def\Top{\bold T\bold o \bold p}
\def\wTop{\text{\rm w}\bold T}
\def\wT{\text{\rm w}\bold T}
\def\Sp{\bold S\bold p}
\def\vo{\varOmega}
\def\vs{\varSigma}
\def\smsh{\wedge}
\def\flush{\flushpar}
\def\id{\text{id}}
\def\dbslash{/\!\! /}
\def\codim{\text{\rm codim\,}}
\def\:{\colon}
\def\holim{\text{\rm holim\,}}
\def\hocolim{\text{\rm hocolim\,}}
\def\hodim{\text{\rm hodim\,}}
\def\hocodim{\text{hocodim\,}}
\def\Bbb{\mathbb}
\def\bold{\mathbf}
\def\Aut{\text{\rm Aut}}
\def\cal{\mathcal}
\def\Sec{\text{\rm sec}}
\def\Secst{\text{\rm sec}^{\text{\rm st}}}
\def\maps{\text{\rm map}}
\setcounter{tocdepth}{1}
\tableofcontents
\addcontentsline{file}{sec_unit}{entry}

\section{Introduction \label{intro}}
In this article an intersection theory
on manifolds is developed using the techniques of algebraic topology. 
The ``cycles'' in our theory will be maps between manifolds,
whereas ``intersections'' will live in the
homotopy groups of a certain Thom space. In order to give the
obstructions a geometric interpretation, one must
identify the homotopy of this Thom space with a suitable
bordism group. Making this identification requires transversality.
Consequently, if one is willing to forgo a geometric
intepretation and work exclusively with Thom spaces, 
transversality can be dispensed with altogether. 
We emphasize this point  for the following reason:
although we will not pursue the matter here, 
our methods straightforwardly extend to give an intersection theory
for Poincar\'e duality spaces, even though the usual transversality 
results are known to fail in this wider context.
\medskip

Suppose
$N$ is an $n$-dimensional compact manifold equipped with
a closed submanifold $$Q \subset N$$ of dimension $q$. 
Given a map $f\:P \to N$, where $P$ is a closed manifold of 
dimension $p$, we ask
for necessary and sufficient conditions insuring that $f$ is homotopic
to a map $g$ whose image is disjoint from $Q$. We call these data
an {\it intersection problem.} The situation is depicted by the diagram
$$
\SelectTips{cm}{}
\xymatrix{
& N - Q \ar[d] \\
P \ar[r]_f\ar@{..>}[ur]
& N \, 
}
$$
where we seek to fill in the dotted arrow by 
a map making the diagram homotopy commute.
Note that transversality implies the above problem always has a solution 
when $p + q < n$. 
\medskip

When $f$ happens to be an embedding (or immersion), one typically requires 
a deformation of $f$ through isotopies (resp.\ regular homotopies). 
This version of the problem was studied
by Hatcher and Quinn \cite{H-Q}, who
approached it geometrically using the methods of bordism theory
(see also the related papers by 
Dax \cite{Dax}, Laudenbach \cite{Laud} and Salomonsen \cite{Salo}).
Since many of the key proofs in \cite{H-Q} are just  sketched, we feel it
is useful to give independent homotopy theoretic proofs of their
results.  Also, some of our steps, such as the {\it Complement Formula}
in Section 5, should be of independent interest. 
\medskip

We now summarize our approach. Let $E \to P$ be the fibration
given by converting the inclusion $N-Q\to N$ into a fibration and then
pulling the latter back along $P$. Then the desired lift 
exists if and only if $E \to P$ admits a section.
Step one is to produce an obstruction whose vanishing
guarantees the existence of such a section. 

In a certain range of dimensions,
it turns out that the complete obstruction to finding
a section has been known
for at least 35 years: it goes by the name of
{\it stable cohomotopy Euler class}
(see e.g. Crabb \cite[Ch.\ 2]{Crabb}, who attributes the ideas
to various people, notably  Becker 
\cite{Becker}, \cite{Becker2} and Larmore \cite{Larmore}).

The second step is to equate the stable cohomotopy Euler class,
a `cohomological' invariant, with a {\it stable homotopy Euler characteristic},
a `homological' one. This is achieved using a version of Poincar\'e 
duality which appeared in \cite{Klein_dualizing}.
The characteristic lives in a homotopy group
of a certain spectrum.  

The third step of the program is to  identify
the spectrum in step two as a Thom spectrum. 
The idea here, which we believe is new, 
is to give an explicit 
homotopy theoretic model
for the complement of the inclusion $Q \subset N$ in a certain
stable range. We exhibit this model in Theorem \ref{complement}
(this is the ``Complement Formula''
alluded to above).

The final step, which is optional, 
is to relate the Thom spectrum of step three with
a twisted bordism theory (cf.\ the next paragraph). 
This is a standard application
of the Thom transversality theorem. 
As pointed out above, this step is omitted in
the case of an intersection problem involving
Poincar\'e duality spaces.
This completes our outline of the program.
We now proceed to state our main results in the manifold
case. This will require some preparation. 

We first
make some well-known remarks about the relationship between 
bordism and Thom spectra. Suppose $X$ is
a space equipped with  a vector bundle $\xi$ of rank $k$.
Consider triples
$$
(M,g,\phi)\, ,
$$
in which $M$ is a closed smooth manifold
of dimension $n$ equipped with map $g\: M \to X$ 
and $\phi$ is stable vector bundle
isomorphism between the normal bundle of $M$ and the pullback $g^*\xi$.
The set of equivalence classes of these under the
relation of bordism defines an abelian group 
$$
\Omega_n(X;\xi)\, ,
$$
in which addition is given by disjoint union. 
This is the {\it bordism group of $X$ with coefficients 
in $\xi$} in degree $n$.

The {\it Thom space} $T(\xi)$
is the one point compactification of the total space of $\xi$.
If $\epsilon^j$ denotes the rank $j$ trivial bundle over $X$, then
there is an evident map $$\Sigma T(\xi\oplus \epsilon^j) \to 
T(\xi\oplus \epsilon^{j+1})$$
which gives the collection $\{T(\xi\oplus \epsilon^j )\}_{j\ge 0}$ 
the structure of a (pre-)spectrum. Its associated $\Omega$-spectrum
is called the {\it Thom spectrum} of $\xi$,
which we denote by $$X^\xi\, .$$

The Thom-Pontryagin construction defines a homomorphism 
$$
\Omega_n(X;\xi) \to \pi_{n+k}(X^\xi)
$$
which is an isomorphism by transversality. 
These remarks apply equally as well in the more
general case when $\xi$ is a virtual vector bundle of rank $k$.

More generally, if $\xi$ is a stable spherical fibration, 
the Thom-Pontryagin homomorphism
is still defined, but can fail to be an isomorphism. The deviation
from it being an isomorphism is detected by the surgery theory 
$L$-groups of $(\pi_1(X),w_1(\xi))$ (cf.\ Levitt \cite{Levitt}, 
Quinn \cite{Quinn}, Jones \cite{Jones}, Hausmann and Vogell \cite{H-V}).
In essence, the $L$-groups detect the 
failure of transversality.
\medskip

We now return to our intersection problem.
Let $i_Q\: Q\subset N$ denote the inclusion. 
Define
$$
E(f,i_Q)
$$
to be {\it homotopy fiber product} of 
$f$ and $i_Q$ (also known as the {\it homotopy pullback}). 
Explicitly, a point of $E(f,i_Q)$ is a triple
$$(x,\lambda,y)$$ 
in which $x\in P$, $y \in Q$ and $\lambda\: [0,1]\to N$
is a path such that $\lambda(0) = f(x)$ and $\lambda(1) = y$.

Define a virtual vector bundle $\xi$ over $E(f,i_Q)$ as follows: 
let 
$$
j_P \: E(f,i_Q) \to P \qquad \text{ and } \qquad j_Q\: E(f,i_Q) \to Q
$$
be the 
forgetful maps and let 
$$
j_N \: E(f,i_Q) \to N
$$ 
be the map
given by $(p,\lambda,q) \mapsto \lambda(1/2)$. Then $\xi$ is
defined as the  rank $n-p-q$  virtual bundle 
$$
(j_N)^*\tau_N - (j_P)^*\tau_P - (j_Q)^*\tau_Q \, ,
$$
where, for example, $\tau_N$ denotes
the tangent bundle of $N$ and $(j_N)^*\tau_N$
is its pullback along $j_N$.

\begin{bigthm} \label{main}
Given an intersection problem,
there is an  obstruction
$$
\chi(f,i_Q) \in \Omega_{p + q - n}(E(f,i_Q); \xi)\, ,
$$
which vanishes when $f$ is homotopic to a map with
image disjoint from $Q$.

Conversely, if $p + 2q+ 3 \le 2n$
and $\chi(f,i_Q)=0$, then the intersection
problem has a solution: there is a homotopy from $f$ to a map
with image disjoint from $Q$.
\end{bigthm}

Our second main result identifies the homotopy
fibers of the map of mapping spaces
$$
\text{\rm map}(P,N-Q) \to \text{\rm map}(P,N)
$$
in a range.
Choose a basepoint $f\in  \text{\rm map}(P,N)$. 
Then we have a homotopy fiber sequence
$$
{\cal F}_f \to \text{\rm map}(P,N-Q) \to \text{\rm map}(P,N)\, ,
$$
where ${\cal F}_f$ denotes the homotopy fiber at $f$.

\begin{bigthm} \label{families}
Assume $\chi(f,i_Q)$ is trivial. Then
there is a  $(2n{-}2q{-}p{-}3)$-connected map
$$
{\cal F_f} \to \Omega^{\infty{+}1} E(f,i_Q)^{\xi} \, ,
$$
where the target is the
loop space of the zero${}^\text{th}$ space of the Thom spectrum
$E(f,i_Q)^{\xi}$. 
\end{bigthm}

Theorems \ref{main} and \ref{families} are the main
results of this work. In \S9-11, we give applications
of these results to fixed point theory, embedding theory 
and linking problems.

\subsection*{Families} Theorem \ref{families} yields
an obtruction to removing intersections in families.
Let
$$
F\:P \times D^j \to N
$$
be a $j$-parameter family of maps whose restriction to
$P \times S^{j-1}$ has disjoint image with $Q$, and
whose restriction to $P \times *$ is denoted by $f$, where
$*\in S^{j-1}$ is the basepoint. Assume $j >0$.

The adjoint of $F$ determines a based map of pairs
$$
(D^j,S^{j-1}) \to  (\text{\rm map}(P,N), \text{\rm map}(P,N-Q))
$$
whose associated homotopy class gives rise to  
an element of $\pi_{j-1}({\cal F}_f)$.  By Theorem \ref{families},
this class is determined by its 
image in the abelian group 
$$
\pi_{j-1}(\Omega^{\infty+1}E(f,i_Q)^\xi)\,\, \cong \,\,
\Omega_{j+p+q-n}(E(f,i_Q);\xi)\, $$
provided $j + p + 2q + 3 \le 2n$. Denote this element by 
 $\chi^{\text{fam}}(F)$. 

\begin{bigcor} Assume $0 < j \le 2n-2q-p-3$. Then the 
family $F \: P\times D^j \to N$ is homotopic rel $P \times S^{j-1}$ 
to a family having disjoint image with $Q$ if and only if
$\chi^{\text{\rm fam}}(F)$ is trivial.
\end{bigcor}

\subsection*{Additional remarks}
We will not deal with self intersection problems
here as that was in effect handled by the first author
in \cite{Klein_compression}.

Some of the machinery developed below has recently been applied
by M.\ Aouina \cite{Aouina} to  identify the homotopy type of
the moduli space of thickenings of a finite complex
in the metastable range.

We plan two sequels to this paper. One
 will consider a multirelative version
of the theory, in which we develop and obstruction theory for
deforming maps off of more than one submanifold. The other
sequel will develop tools to study periodic points of
dynamical systems.
\medskip

\begin{acks} We are much indebted to Bill Dwyer for discussions
that motivated this work. We were also to a great extent inspired
by the ideas of the Hatcher-Quinn paper \cite{H-Q}.
\end{acks}

\begin{out} Section 2 sets forth language. Section 3 contains
various results about section spaces, and we
define the stable cohomotopy Euler class. Most of this material is
probably classical. In 
Section 4 we use a version of Poincar\'e duality to define the
stable homotopy Euler class. Section 5 contains
a ``Complement Formula'' which identifies the homotopy type
of the complement of a submanifold in a stable range.
Section 6 contains the proof of Theorem
\ref{main} and Section 7 the proof of Theorem \ref{families}. Section 8 gives
an alternative definition of the main invariant which doesn't
require $i_Q \: Q \to N$ to be an embedding. Section 9 describes 
a generalized linking invariant based on our intersection invariant.
Section 10 applies our results to fixed point problems. Sections 11-12 show how
our results, in conjunction with a result of Goodwillie and first
author, can be used to deduce the intersection theory of Hatcher and Quinn.
\end{out}

\section{Language}

\subsection*{Spaces}
All spaces will be compactly generated, and
products are to be re-topologized using the compactly
generated topology. Mapping spaces are to be given
the compactly generated, compact open topology.
A {\it weak equivalence} of
spaces denotes a (chain of) weak homotopy equivalence(s).

Some connectivity conventions: 
the empty space is $(-2)$-connected.
Every nonempty space is $(-1)$-connected. A nonempty
space $X$ is $r$-connected for $r \ge 0$ if $\pi_j(X,*)$ 
is trivial for $j \le r$ for all choices of basepoint $* \in X$.
A map $X \to Y$ of nonempty spaces is  $(-1)$-connected
and is $0$-connected if it is surjective on path components.
It is $r$-connected for $r > 0$ if all of its homotopy fibers 
are $(r-1)$-connected.

When speaking of manifolds, we work exclusively in the
smooth ($C^\infty$) category. However, all results
of the paper hold equally well in the PL and topological categories.

\subsection*{Spaces as classifying spaces}
Suppose that $Y$ is a connected based space. 
The simplicial total singular complex of $Y$ is a based simplicial
set. Take its Kan simplicial loop group. Define
$G_Y$ to be its geometric realization. Then $G_Y$ is
a topological group and there is
a functorial weak equivalence 
$$Y \simeq BG_Y \, ,$$
where $BG_Y$ is the classifying space of $G_Y$.
To emphasize
that $G_Y$ is a group model for the loop space, we usually abuse
notation and rename 
$$
\Omega Y:= G_Y\, .
$$
Thus, $Y$ is identified with $B\Omega Y$.

\subsection*{The thick fiber of a map}
Let  
$$
f\:X\to B
$$
be a map, where $B = BG$ is the classifying space
of a topological group $G$ which is the realization
of a simplicial group.

The 
{\it thick fiber} of $f$ is the space 
$$
F := \text{pullback}(EG \to B \leftarrow X)
$$
where $EG \to B$ is a universal principal $G$-bundle.
Let $G$ act on the product $EG \times X$ by means of
its action on the first factor. This action leaves 
the subspace $F \subset EG \times X$ setwise invariant,
so $F$ comes equipped with a $G$-action. 

Furthermore, as $EG$ is contractible, 
$F$ has the homotopy type of the homotopy fiber of $f$.
Observe that $X$ has the weak homotopy type of
the Borel construction $EG\times_G F$ and  
$f\:X\to B$ is then identified up to homotopy 
as the fibration $EG \times_{G} F \to B$.

\subsection*{Naive equivariant spectra} We will be using a low tech
version of equivariant spectra, which are defined
over any topological group.

Let $G$ be as above. A {\it (naive) $G$-spectrum} $E$ 
consists of based
(left) $G$-spaces $E_i$ for $i \ge 0$, and  
equivariant based maps $\Sigma E_i \to E_{i{+}1}$
(where we let $G$ act  trivially on the suspension coordinate of $\Sigma E_i$).
A {\it morphism} $E \to E'$ of $G$-spectra consists of
maps of based spaces $E_i \to E'_i$ which are compatible with
the structure maps.  A {\it weak equivalence}  of
$G$-spectra is a map inducing an isomorphism on
homotopy groups. $E$ is an {\it $\Omega$-spectrum} if
the adjoint maps $E_i \to \Omega E_{i{+}1}$ are weak
equivalences.
We will for the most part
assume that our spectra are $\Omega$-spectra. 
If $E$ isn't an $\Omega$-spectra,
we can functorially approximate it by one: $E \overset \sim\to E'$, where
$E'_i$  is the homotopy colimit of the diagram of
$G$-spaces $\{\Sigma^kE_{i+k}\}_{k\ge 0}$.  We use the notation
$\Omega^\infty E$ for $E'_0$, and by slight abuse of language,
we call it the {\it zero${}^{\rm th}$ space} of $E$.

If $X$ is a based $G$-space, then its 
{\it suspension spectrum} $\Sigma^{\infty} X$ is a $G$-spectrum with
$j$-th space $Q(S^j \smsh X)$,
where $Q = \Omega^{\infty}\Sigma^{\infty}$ is the stable
homotopy functor.

The {\it homotopy orbit spectrum} $$E_{hG}$$ of $G$ acting on $E$ is
the spectrum whose $j^{\text{\rm th}}$ space is
the orbit space of $G$ acting diagonally 
on the smash product $E_j \smsh EG_+$. The structure maps
in this case are evident.

\section{The stable cohomotopy Euler class}

Suppose $p\:E \to B$ is a Hurewicz fibration 
over a connected space $B$.
We seek a generalized cohomology theoretic 
obstruction to finding a section.

The {\it fiberwise suspension} of $E$ over $B$ 
is the double mapping cylinder
$$
S_B E  :=  B\times 0  \cup E \times [0,1] \cup B \times 1\, .
$$
This comes equipped with a map $S_Bp\: S_B E \to B$ which
is also a fibration (cf.\ \cite{Strom}). 
If $F_b$ denotes a fiber of $p$ at $b\in B$, then the fiber
of $S_B p$ at $b$ is
$S F_b$, the unreduced suspension of $F_b$. 

Let $$
s_-,s_+\:B \to S_B E 
$$
denote the sections given by the inclusions of $B\times 0,B\times 1$
into the double mapping cylinder.

\begin{prop} \label{section_lemma} If $p\:E \to B$ admits a section, then
$s_-$ and $s_+$ are vertically homotopic.

Conversely,
assume $p\:E\to B$ is  $(r{+}1)$-connected and
$B$ homotopically a retract of a cell complex 
with cells in dimensions $\le 2r + 1$.
If  $s_-$ and $s_+$ are vertically
homotopic, then $p$ has a section.
\end{prop}
\bigskip

\begin{proof} Assume $p\:E\to B$ has a section $s\: B\to E$.
Apply the functor $S_B$ to $s$ to get a map
$$
S_Bs \: S_B B \to S_B E
$$
and note $S_B B = B \times [0,1]$. Then $S_B s$ is
a vertical homotopy from $s_-$ to $s_+$.

To prove the converse, consider the square
$$
\SelectTips{cm}{}
\xymatrix{
E \ar[r]^p \ar[d]_p & B\ar[d]^{s_+}\\
B \ar[r]_{s_-} & S_B E
}
$$
which is commutative up to preferred homotopy. The square
is homotopy cocartesian. Since the maps out of $E$ are $(r{+}1)$-connected,
we infer via the Blakers-Massey theorem that the square is 
$(2r{+}1)$-cartesian.

Let $\cal P$ be the homotopy pullback of $s_-$ and $s_+$. Then the map 
$E\to {\cal P}$ is $(2r{+}1)$-connected. Furthermore,
a choice of vertical homotopy from $s_-$ to $s_+$ yields  a map 
$B\to {\cal P}$. Using the dimensional constraints on $B$,
we can find a map $s\:B \to E$ which factorizes
$B \to {\cal P}$ up to homotopy. Then $ps$
is homotopic to the identity.
The homotopy lifting property then enables us to deform $s$ to an actual
section of $p$.
\end{proof}

\subsection*{Section spaces}
For a fibration $E\to B$, let
$$
\Sec(E \to B)
$$
denote its space of sections. 
Proposition \ref{section_lemma} gives
criteria for deciding when this space is non-empty. Another way
to formulate it is to consider $\Sec(S_B E\to B)$ as a
space with basepoint $s_-$. Then  
the obstruction of \ref{section_lemma}
is given by asking whether the homotopy class
$$
[s_+] \in \pi_0(\Sec(S_B E \to B))
$$
is that of the basepoint.

\subsection*{Stabilization}
As remarked above, $s_-$ equips the section space 
$$
\Sec(S_B E \to B)
$$ 
with a basepoint, and the fibers $SF_b$ of
$S_B E \to B$ are based spaces
(with basepoint given by the south pole). 
\bigskip

Let $$Q_BS_B E  \to B$$ 
be the effect of applying the stable homotopy
functor $Q = \Omega^\infty \Sigma^\infty$ to each fiber $SF_b$ 
of $S_B E \to B$.

\begin{lem} Assume  $p\: E\to B$ is $(r{+}1)$-connected and 
 $B$ homotopically the retract  of a cell complex
with cells in dimensions $\le k$.
Then the evident map
$$
\Sec(S_B E \to B) \,\, \to \,\, \Sec(Q_BS_B E \to B)
$$
is $(2r {-} k {+}3)$-connected. 
\end{lem}

\begin{proof} For each $b\in B$, the space $SF_b$ is $(r{+}1)$-connected.
The Freudenthal suspension theorem implies that
the  map $SF_b \to QSF_b$ is $(2r{+}3)$-connected. From this we infer
that the map $S_B E \to Q_BS_B E$ is $(2r{+}3)$-connected. The
result now follows from elementary obstruction theory.
\end{proof}

We call $ \sec(Q_BS_B E \to B)$ the {\it stable section
space} of $S_B E \to B$ and we change its notation to
$$
\Secst(S_B E \to B)\, .
$$
In fact, the stable section space is
the zero${}^{\rm th}$ space of a spectrum whose $j$-th space
is the stable section space of a fibration $E_j \to B$
in which the fiber at $b\in B$ is $Q\Sigma^j SF_b$. In particular,
the set of path components of $\Secst (S_B E \to B)$ has the structure of
an abelian group.
In what follows, we regard $s_-,s_+$ as points
of $\Secst (S_B E \to B)$.

\begin{defn} The {\it stable cohomotopy Euler class} of 
$p\:E\to B$ is given by
$$
e(p) := [s_+] \in \pi_0(\Secst (S_B E \to B))\, .
$$
\end{defn}

\begin{cor} If $p$ has a section, then $e(p)$ is trivial.
Conversely, assume $p\:E\to B$ is $(r{+}1)$-connected and
$B$ is homotopically the retract of a cell complex with cells in dimensions
$\le 2r{+}1$. If $e(p) = 0$, then 
$p$ has a section.
\end{cor}

We will also require an alternative description of the
homotopy type of $\sec(E\to B)$ in a range.
 For a space $X$ equipped with two points $-,+\in X$, let
$$
\Omega^\pm X
$$
be the space of paths $\lambda\:[0,1] \to X$ such that $\lambda(0) = -$
and $\lambda(1) = +$. When $X = SY$, and $\pm$ are the poles of $SY$,
we obtain a natural map
$$
Y \to \Omega^\pm QS Y
$$
which maps a point $y$ to the path $[0,1] \to QSY$ by 
$t\mapsto t\smsh y$, where $t\smsh y \in SY$ is considered as
a point of $QSY$ in the evident way.

Next, suppose that $E\to B$ is a fibration. Then the 
associated fibration
$$
Q_B S_B E \to B
$$
has $QSF_b$ as its fiber at $b\in B$. So we have a map
$$
F_b \to \Omega^\pm QSF_b\, .
$$
Let 
$$
\Omega^\pm_B Q_B S_B E \to B
$$
be the fibration whose fiber at $b$ is the space
$\Omega^\pm QSF_b$. Then the above yields a map of
section spaces
$$
\Sec (E \to B) \to \Sec (\Omega^\pm_B Q_B S_B E \to B)\, .
$$

\begin{lem} \label{section_type_1}  Assume $E\to B$ is $(r{+}1)$-connected
and $B$ is homotopically the retract of a cell complex with 
cells in dimensions $\le k$. Then the map 
$$
\Sec (E \to B) \to \Sec (\Omega^\pm_B Q_B S_B E \to B)
$$
is $(2r{+}1 - k)$-connected.
\end{lem}

\begin{proof} For each $b\in B$, the map of fibers
$$
F_b \to \Omega^\pm Q S F_b
$$
factors as a composite
$$
F_b \to \Omega^\pm SF_b \to \Omega^\pm QSF_b\, .
$$
The first map in the composite is $(2r{+}1)$-connected
(by the Blakers-Massey theorem) and the second map is
is $(2r{+}2)$-connected (by Freudenthal's suspension theorem).
Hence, the composed map is $(2r{+}1)$-connected. We infer by the
five lemma that the map $E \to \Omega^\pm_B Q_B S_B E$ is also
$(2r{+}1)$-connected. Taking section spaces then reduces 
the connectivity by $k$.
\end{proof}

\begin{lem}\label{section_type_2} Fix a section 
$s$ of the fibration $E \to B$.
Then with respect to this choice, there is a preferred weak equivalence
$$
\Sec (\Omega^\pm_B Q_B S_B E \to B) \,\, \simeq \,\, 
\Omega \Sec (Q_B S_B E \to B)\, .
$$
\end{lem}

\begin{proof} The fiber of $(\Omega^\pm_B Q_B S_B E \to B)$ at $b\in B$
is the space $\Omega^\pm QSF_b$. The hypothesis that
$E\to B$ is equipped with a section 
shows that $F_b$ is based and therefore $\Omega^\pm QSF_b$
is also based using the map $F_b \to \Omega^\pm QSF_b$.

A point of $\Omega^\pm QSF_b$ 
is a path in $QSF_b$ having fixed endpoints.
Given another point of $\Omega^\pm QSF_b$, we get
another path having the same
endpoints. Now form the loop which starts by traversing the first
path and returns by means of the second path. So we get a map
$$
\Omega^\pm QSF_b \to \Omega QSF_b 
$$
which is a weak equivalence (an inverse
weak equivalence is given by mapping a loop in $QSF_b$ to
the path given by concatenating the base path with the given loop).
This weak equivalence then induces
a weak equivalence 
$$
\Sec (\Omega^\pm_BQ_B S_BE \to B) \to \Sec (\Omega_BQ_B S_B E \to B) 
$$
where $\Omega_B$ is the fiberwise loop space functor.
Now use the evident homeomorphism
$$
\Sec (\Omega_BQ_B S_B E\to B) \cong  \Omega \Secst(S_BE \to B) 
$$
to complete the proof.
\end{proof}

Assembling the above lemmas, we conclude

\begin{cor}\label{section_type_3}
Let $E\to B$ be a fibration equipped with section. 
Assume $E \to B$ is $(r{+}1)$-connected and that
$B$ is homotopically the retract of a cell complex 
whose cells have dimension $\le k$.
Then the map
$$
\Sec(E\to B) \to \Omega \Secst (S_BE\to B)
$$
is $(2r+1-k)$-connected.
\end{cor}

\section{The stable homotopy Euler characteristic}

Let $B$ be a connected based space. We identify 
$B$ with the classifying space of the topological
group $\Omega B$ described in $\S2$.
Let $E \to B$ be a fibration, and let $F$ be its thick fiber (\S2).
Take its unreduced suspension $SF$. Then $SF$ is a based $\Omega B$-space.

Assume now that $B$ is a closed manifold of dimension $d$.
Let $\tau_B$ be its tangent bundle, and let $S(\tau_B +\epsilon)$ be
the fiberwise one point compactification of $\tau_B$. 
Define $$S^{\tau_B}$$ 
to be its thick fiber. This is a based $\Omega B$-space.

Define $$S^{-\tau_B}$$ to be the mapping spectrum
$\text{map}(S^\tau_B,S^0)$, i.e., the spectrum whose
$j^{\text{\rm th}}$ space consists of the stable based maps from $S^{\tau_B}$
to $S^j$. Then $S^{\tau_B}$ is an $\Omega B$-spectrum whose
underlying unequivariant homotopy type is that of a $(-d)$-sphere.

Give the smash product $S^{-\tau_B} \smsh SF$ 
the diagonal action of
$\Omega B$. 
Let $S^{-\tau_B}\smsh_{h\Omega B} SF$ be its homotopy orbit spectrum.

\begin{thm}[``Poincar\'e Duality''] \label{duality} 
There is a preferred weak equivalence of
infinite loop spaces
$$
\Omega^\infty (S^{-\tau_B}\smsh_{h\Omega B} SF) \,\, \simeq \,\, 
\Secst(S_B E \to B)\, .
$$
In particular, there is an preferred isomorphism of abelian groups
$$
\pi_0(S^{-\tau_B}\smsh_{h\Omega B} SF) \,\, \cong \,\,  \pi_0(\Secst(S_B E \to B))\, .
$$
\end{thm}

\begin{rem} A form of this statement which resembles classical
Poincar\'e duality is given
by thinking of the right side as ``cohomology'' with coefficients
in the ``cosheaf'' of spectra ${\cal E}$
over $B$ whose stalk at $b\in B$ is
the spectrum $\Sigma^\infty (S F_b)$. 
Then symbolically, the result identifies cohomology
with twisted homology:
$$
H_\bullet(B;S^{-\tau_B} \otimes {\cal E}) \simeq  H^{\bullet}(B;{\cal E})\, ,
$$
where we interpret the displayed tensor product as fiberwise smash product.
\end{rem}

\begin{proof}[Proof of Theorem \ref{duality}]
The theorem is actually a special case of the main
results of \cite{Klein_dualizing}. There,
for any $\Omega B$-spectrum $W$, we constructed
a weak natural transformation 
$$
S^{-\tau_B} \smsh_{h\Omega B} W \to W^{h\Omega B}
$$
called the {\it norm map,} which was subsequently shown
to be a weak equivalence for every $W$ (cf.\ th.\ D and 
cor.\ 5.1 of \cite{Klein_dualizing}). The target of the norm map is
the homotopy fixed points of $\Omega B$ acting on $W$.
Recall that when $W$ is an $\Omega$-spectrum, $W^{h\Omega B}$ is the spectrum
whose $j^{\text{\rm th}}$ space is the section space of the fibration
$E\Omega B \times_{\Omega B} W_j \to B$.

Specializing the norm map to the $\Omega B$-spectrum 
$W = \Sigma^\infty SF$, source of the norm  map is identified with
$S^{-\tau_B} \smsh_{h\Omega B} SF$ whereas its target is 
identified with the spectrum whose associated infinite loop space
is $\Secst(S_B E \to B)$.
\end{proof}

\begin{defn} The {\it stable homotopy Euler characteristic} of
$p\:E \to B$ is the class
$$
\chi(p) \in \pi_0(S^{-\tau_B}\smsh_{h\Omega B} SF)
$$
which corresponds to the stable cohomotopy Euler class
$e(p)$ via the isomorphism of  Theorem \ref{duality}.
That is, $\chi(p)$ is the Poincar\'e dual of $e(p)$.
\end{defn}

\begin{cor} \label{chi-cor} Assume $p\: E\to B$ is $(r{+}1)$-connected, 
$B$ is a closed manifold dimension $d$ and
$d  \le 2r+1$. Then $p$ has a section if and only if 
$\chi(p)$ is trivial.
\end{cor}

\section{The complement formula}

Suppose
$$
i_Q\:Q^q \subset N^n
$$
is the inclusion of a closed connected submanifold. We will 
also assume that $N$ is connected.

Choose a basepoint $* \in Q$. Then $N$ gets a basepoint.
Fix once and for all identifications
$$
Q \simeq B\Omega Q \qquad  N \simeq B\Omega N\, .
$$
We also have a homomorphism $\Omega Q \to \Omega N$, such
that application of the classifying space functor 
yields $i_Q\: Q\to N$ up to homotopy.

Let
$$
F := \text{thick fiber}(N-Q \to N)
$$
be the thick fiber of the inclusion $N-Q \to N$ taken
at the basepoint. Then $F$ is a $\Omega N$-space whose unreduced
suspension $SF$ is a based $\Omega N$-space. Consequently, the suspension
spectrum $\Sigma^\infty SF$ has the structure of an $\Omega N$-spectrum.
We will identify the equivariant homotopy
type of this spectrum.

Let $\nu_{Q\subset N}$ denote the normal bundle of $i_Q\:Q\to N$. 
Form the fiberwise one-point compactification of this vector bundle
to obtain a sphere bundle over $Q$ equipped with a preferred section
at $\infty$. Let
$$
S^{\nu_{Q\subset N}}
$$
denote its thick fiber. This is, up to homotopy, a sphere whose 
dimension coincides with the rank of $\nu_{Q\subset N}$. 
Then $S^{\nu_{Q\subset N}}$ comes equipped with an $\Omega Q$-action.

Another way to construct $S^{\nu_{Q\subset N}}$ which emphasizes its
dependence only  on the homotopy class of $i_Q$ is as follows: let
$S^{\tau_N}$ be the thick fiber of the fiberwise one point compactification
of the tangent bundle of $N$. This is an $\Omega N$-spectrum, and
therefore and $\Omega Q$-spectrum by restricting the action.

Similarly, using the tangent bundle $\tau_Q$ of $Q$, we obtain
a based space $S^{\tau_Q}$ equipped with an $\Omega Q$-action.
Let
$$
S^{\tau_N - \tau_Q} 
$$
be the spectrum whose $j^{\text{\rm th}}$ space
consists of the stable 
based maps from $S^{\tau_Q}$ to the $j$-fold reduced suspension
of $S^{\tau_N}$. 
Then $S^{\tau_N - \tau_Q}$ is a an $\Omega Q$-spectrum.

Using the stable bundle isomorphism
$$
\nu_{Q\subset N} \cong i_Q^*\tau_N - \tau_Q
$$
it is elementary to check that the suspension spectrum 
of $S^{\nu_{Q\subset N}}$ (an $\Omega Q$-spectrum) has the same
weak equivariant homotopy type as  $S^{\tau_N - \tau_Q}$.

\begin{thm}[Complement Formula]\label{complement} 
There is a preferred equivariant
weak equivalence of $\Omega N$-spectra
 $$
\Sigma^\infty SF \,\, \simeq \,\, 
S^{\tau_N -\tau_Q} \smsh_{h \Omega Q} (\Omega N)_+\, ,
$$
where the right side is the homotopy orbits of $\Omega Q$
acting diagonally on the smash product
$S^{\tau_N -\tau_Q} \smsh (\Omega N)_+$.
Alternatively, it 
is the effect of inducing
$S^{\tau_N - \tau_Q}$ along the homomorphism
$\Omega Q \to \Omega N$ in a homotopy invariant way.
\end{thm}

\begin{rem} This result recovers the homotopy type of 
$S_N(N-Q)$ as a space over $N$ in the stable range. Namely,
the Borel construction applied to 
$S^{\tau_N -\tau_Q} \smsh_{h \Omega Q} (\Omega N)_+$ gives
a family of spectra over $N$, and the result says that
this family coincides up to homotopy with the fiberwise
suspension spectrum of $S_N(N-Q)$ over $N$.
\end{rem}

\begin{proof}[Proof of \ref{complement}] Let $\nu := \nu_{Q\subset N}$ 
denote the normal bundle of $Q$ in $N$.
Using a choice of tubular neighborhood, we have a 
homotopy cocartesian square
$$
\SelectTips{cm}{}
\xymatrix{
S(\nu) \ar[r]\ar[d] & N - Q\ar[d] \\
D(\nu) \ar[r] & N \, .
}
$$
There is then a weak equivalence of homotopy colimits
\begin{equation} \label{N-pushout}
\begin{CD} 
\hocolim (D(\nu) \leftarrow S(\nu) \to N) 
\overset \sim \to
\hocolim (N \leftarrow N-Q \to N)\, .
\end{CD}
\end{equation}
Each space appearing in \eqref{N-pushout} is a space over $N$. 
Take the thick fiber over $N$ of each of these spaces to get
an equivariant weak equivalence
\begin{equation} \label{fiber_pushout}
\hocolim (\tilde D(\nu) \to \tilde S(\nu) \leftarrow *) 
\overset\sim \to 
\hocolim (* \to F \leftarrow *) =: S F\, ,
\end{equation}
where $F$ is the thick fiber of
$N-Q \to N$. The proof will be completed by identifying
the domain of \eqref{fiber_pushout}.

If $\tilde Q$ denotes the thick fiber of $i\:Q \to N$,
then the domain of \eqref{fiber_pushout}
is, by definition, the Thom space of the
pullback of $\nu$ along the map $\tilde Q \to Q$.

We can also identify $\tilde Q$ with the homotopy orbits
$\Omega Q$ acting on $\Omega N$:
$$
\tilde Q \simeq E\Omega Q \times_{\Omega Q}\Omega N\, .
$$
This space comes equipped with a spherical
fibration given by 
\begin{equation}\label{fibration}
E\Omega Q \times_{\Omega Q}(\Omega N \times S^\nu) \to 
E\Omega Q \times_{\Omega Q}(\Omega N \times *) = \tilde Q
\end{equation}
where $S^\nu$ denotes the thick fiber of the spherical fibration
given by fiberwise one point compactifying $\nu$. The spherical
fibration \eqref{fibration} comes equipped with a preferred section
(coming from the basepoint of $S^\nu$.
It is straightforward to check that this fibration coincides with 
the fiberwise one point compactification of the pullback of $\nu$ to 
$\tilde Q$.

Consequently, the Thom space of the pullback of $\nu$ to 
$\tilde Q$ coincides up to homotopy with
the effect
collapsing the preferred section of \eqref{fibration}
to a point. But the effect of this collapse 
this yields
$$
(\Omega N)_+ \smsh_{h\Omega Q} S^\nu \, .
$$
Hence, what we've exhibited is an $\Omega N$-equivariant 
weak equivalence of based spaces
$$
SF \,\, \simeq \,\, (\Omega N)_+ \smsh_{h\Omega Q} S^\nu \, .
$$
The proof is completed by taking the suspension spectra
of both sides and recalling that $\Sigma^\infty S^\nu$ 
is  $S^{\tau_N - \tau_Q}$.
\end{proof}

\section{Proof of Theorem \ref{main}}

Returning the the situation of the introduction, suppose 
$$
\SelectTips{cm}{}
\xymatrix{
& N - Q \ar[d] \\
P \ar[r]_f\ar@{..>}[ur]
& N 
}
$$
is an itersection problem. 
As already mentioned, the obstructions to lifting
$f$ up to homotopy coincide with the obstructions to sectioning
the fibration
$$
p\:E \to P
$$
where $E$ is the homotopy fiber product of $P \to N \leftarrow N-Q$.

Choose a basepoint for $P$. Then $N$ gets a basepoint via $f$.
The the thick fiber $p$ at the basepoint is identified
with the thick fiber of $N-Q \to N$. Call the thick fiber
of the latter $F$. Then $F$ is an $\Omega N$-space; using the
homomorphism $\Omega P \to \Omega N$, we see that
$F$ is also an $\Omega P$-space.

\begin{lem}\label{N-Q->N}
The map $N-Q \to N$ is $(n-q-1)$-connected.
\end{lem}

\begin{proof} Using the tubular
neighborhood theorem, $N-Q \to N$ is the cobase
change up to homotopy of the spherical fibration
 $S(\nu_{Q\subset N}) \to Q$ of the normal
bundle of $Q$. The fibers of this fibration are
spheresd of dimension $n-q-1$, so the fibration
is that much connected. Then $N-Q\to N$ is also
$(n{-}q{-}1)$-connnected because cobase change
preserves connectivity.
\end{proof}

\begin{cor} \label{conn-E->P} 
The map $E\to P$ is also $(n-q-1)$-connected.
\end{cor}

\begin{proof} $E\to P$ is the base change
of the the $(n-q-1)$-connected map $N-Q\to N$
converted into a fibration. 
The result follows from the fact that base change
preserves connectivity.
\end{proof}

We will now apply Corollary \ref{chi-cor}.
For this we note that the manifold
$P$ is homotopically a cell complex of
dimension $p$, Consequently, if 
$$
p \le 2(n - q - 2) + 1 = 2n - 2q -3\, ,
$$
\ref{chi-cor} implies that a section exists
if and only if the stable homotopy Euler characteristic
$$
\chi(p) \in \pi_0(S^{-\tau_P}\smsh_{h\Omega P} SF)
$$
is trivial.

To complete the proof, 
we will need to identify the homotopy type of the spectrum
\begin{equation}\label{SF-thom}
S^{-\tau_P}\smsh_{h\Omega P} SF\, .
\end{equation}
By the Complement Formula \ref{complement}, 
there is a preferred weak equivalence
of $\Omega N$-spectra
$$
\Sigma^\infty SF \,\, \simeq \,\,  S^{\tau_N -\tau_Q} 
\smsh_{h\Omega Q} (\Omega N)_+\, .
$$
Substituting this identification into \eqref{SF-thom},
 we get a weak equivalence
of spectra
\begin{equation} \label{left-right}
S^{-\tau_P}\smsh_{h\Omega P} SF
\simeq 
S^{-\tau_P} \smsh_{h\Omega P} S^{\tau_N -\tau_Q}  
\smsh_{h\Omega Q} (\Omega N)_+ \, .
\end{equation}

To identify the right side of \eqref{left-right}
as a Thom spectrum, rewrite it again as
$$
S^{\tau_N - \tau_P - \tau_Q} \smsh_{h(\Omega P \times \Omega Q)} (\Omega N)_+\, .
$$
Here, the action of $\Omega P \times \Omega Q$ on 
$$
S^{\tau_N - \tau_P - \tau_Q} = S^{\tau_N} \smsh S^{-\tau_P} \smsh S^{-\tau_Q}
$$ 
is given by having $\Omega P$ act trivially  on $S^{-\tau_Q}$, having
$\Omega Q$ act trivially on $S^{-\tau_P}$, and having
$\Omega P$ and $\Omega Q$ act on $S^{\tau _N}$ 
 by restriction of the $\Omega N$ action.

Clearly, this is the Thom spectrum associated to the
(stable) spherical fibration 
$$
(S^{\tau_N - \tau_P-\tau_Q} \times  E\Omega (P \times Q)) \times_{\Omega P \times \Omega Q} \Omega N
\to  E\Omega (P \times Q)\times_{\Omega (P \times Q)}\Omega N\, .
$$
It is straightforward to check that the base space of this 
fibration 
is weak equivalent to the homotopy fiber product 
$E(f,i_Q)$
described in the introduction. 

Hence, the right side of \eqref{left-right} is just a Thom spectrum of the virtual
bundle $\xi$ over $E(f,i_Q)$ (where $\xi$ is defined as in the introduction). 
Therefore, we get an equivalence of Thom spectra,
$$
S^{-\tau_P}\smsh_{h\Omega P} SF \,\, \simeq \,\, E(f,i_Q)^\xi \, .
$$
Consequently, the stable homotopy Euler characteristic 
becomes identified with an element
of the group
$$
\pi_0(E(f,i_Q)^\xi)
$$
which, by transversality, coincides with the bordism group
$$
\Omega_{p+q-n}(E(f,i_Q);\xi) \, .
$$
Therefore, $\chi(p)$ can be regarded as an element of this bordism group.
The proof of Theorem \ref{main} is then completed by applying  \ref{chi-cor}.

\section{Proof of Theorem \ref{families}}
Recall the weak equivalence
$$
{\cal F}_f \,\, \simeq \,\, \Sec (E\to P)\, ,
$$
where ${\cal F}_f$  is the homotopy fiber of 
$$
 \text{\rm map}(P,N-Q) \to \text{\rm map}(P,N)
$$
at $f\:P\to N$, and 
$$
E \to P
$$
is the fibration in which $E$ is the homotopy pullback
of  
$$
\begin{CD}
P @> f >> N @<\supset << N-Q\, .
\end{CD}
$$
Using \ref{conn-E->P}, we have that $E \to P$ is $(n{-}q{-}1)$-connected.
So by  Corollary \ref{section_type_3},
there is a $(2n-2q-p-3)$-connected map
$$
{\cal F_f} \simeq \Sec(E\to P) \to \Omega \Secst (S_PE\to P)\, .
$$

By \ref{duality}
and the Complement Formula \ref{complement},
there is a 
weak equivalence
$$
\Secst(S_P E \to P)\,\, \simeq \,\,  \Omega^\infty E(f,i_Q)^\xi
$$
where the right side is the Thom spectrum of 
associated with the bundle $\xi$ appearing in the introduction.
Looping this last map, and using the previous identifications, we obtain 
a  $(2n-2q-p-3)$-connected map 
$$
{\cal F}_f \to   \Omega^{\infty+1} E(f,i_Q)^\xi \, .
$$
This completes the proof of Theorem \ref{families}.

\section{A symmetric description}

It is clear from its construction that the stable
homotopy Euler characteristic
depends only upon the homotopy class of $f\: P \to N$
and the isotopy class of $i_Q\: Q \to N$.
Using a different description of the invariant, we will 
explain why it is an invariant of the {\it homotopy class} of $i_Q$.

The new description is more general in that it is
defined not just for inclusions $i_Q\:Q \to N$ but for any map, and
it is symmetric in $P$ and $Q$.

Given maps $f\: P \to N$ and $g\:Q\to N$, consider the intersection problem
$$
\SelectTips{cm}{}
\xymatrix{
& N\times N - \Delta \ar[d] \\
P\times Q \ar[r]_{f\times g}\ar@{..>}[ur]
& N\times N 
}
$$
which asks to find a deformation of $f\times g$ to a map
missing the diagonal $\Delta$ of $N \times N$. Then we
have a stable homotopy Euler characteristic
$$
\chi(f\times g,i_\Delta) \in \Omega_{p+q-n}(E(f\times g,i_\Delta);\xi')
$$
for a suitable virtual bundle $\xi'$ defined as in the introduction. 
A straightforward chasing
of definitions shows there to be a homeomorphism of spaces
$$
E(f\times g,i_\Delta) \,\, \cong \,\, E(f,g) \, . 
$$
Furthermore, if $\xi$ is 
the virtual bundle on $E(f,g)$ defined as in the introduction,
it is clear that $\xi'$ and $\xi$ are equated via this homeomorphism.
The upshot of these remarks is that we can think of $\chi(f\times g,i_\Delta)$
as an element of the bordism group
$$
\Omega_{p+q-n}(E(f,g);\xi)\, .
$$
When $g = i_Q$ is an embedding, this is the same place where
the our originally defined invariant $\chi(f,i_Q)$ lives.

\begin{thm} \label{symmetry}
In fact, the invariants $\chi(f\times i_Q,i_\Delta)$
and $\chi(f,i_Q)$ are equal.
\end{thm}

\begin{rem} Theorem \ref{symmetry} immediately shows that
$\chi(f,i_Q)$ depends only on the homotopy classes of
$f$ and $i_Q$. It also gives extends our intersection invariant 
to the case when $i_Q$ isn't an embedding.
\end{rem}

\begin{proof}[Proof of \ref{symmetry}] (Sketch).
The way to compare the invariants is to consider the commutative
diagram
$$
\SelectTips{cm}{}
\xymatrix{
    & N - Q \ar[r] \ar[d] &  \text{map}(Q,N^{\times 2} - \Delta)\ar[d] \\
P \ar[r]_f\ar@{..>}[ur] \ar@{..>}[rru] & N \ar[r]  & \text{map}(Q,N^{\times 2}) }
$$
where the vertical maps are inclusions and the horizontal
maps are described by mapping a point $x$ to the map
$y \mapsto (x,y)$. The obstruction associated
with short dotted arrow is $\chi(f,i_Q)$, whereas the one associated
with the long one is $\chi(f\times i_Q,i_\Delta)$. The result
now follows from naturality.
\end{proof}

\section{Linking}

Let $P,Q$ and $N$ be compact manifolds where $P$ and $Q$ are
closed.
$$
{\cal D} \subset  \maps(P,N) \times \maps(Q,N)
$$ 
denote the 
{\it discriminant locus} consisting of those pairs of maps $$(f,g)$$ 
such that $f(P) \cap g(Q)$ is non-empty.
Then  $$\maps(P,N) \times \maps(Q,N) - {\cal D}$$ is the space of pairs
$(f,g)$ such that $f$ and $g$ have disjoint images.
\medskip

Consider the commutative diagram
$$
\SelectTips{cm}{}
\xymatrix{
\maps(P,N) \times \maps(Q,N) - {\cal D}
\ar[r] \ar[d]  & \maps(P\times Q, N^{\times 2} - \Delta)\ar[d]  \\
\maps(P,N) \times \maps(Q,N) \ar[r] & 
 \maps(P\times Q, N^{\times 2}) \, ,
}
$$
where the horizontal maps are defined by $(f,g) \mapsto f \times g$,
and the vertical ones are inclusions.

Fix a basepoint $(f,g) \in \maps(P,N) \times \maps(Q,N)$, and define
$$
F_{f,g} := \text{hofiber}(\maps(P,N) \times \maps(Q,N) - {\cal D} \to
\maps(P,N) \times \maps(Q,N))\, ,
$$
where the homotopy fiber is taken at the basepoint.
Likewise, let 
$$
L_{f,g} :=  \text{hofiber}(\maps(P\times Q, (N\times N) - \Delta)
\to  \maps(P\times Q, N \times N))
$$
where the homotopy fiber is taken at the basepoint $f\times g$.
Then we have a map of homotopy fibers
$$
F_{f,g} \to L_{f,g}\, .
$$

According to Theorem \ref{families}
 and the remarks of the previous section, there is a $(2n-p-q-3)$-connected map
$$
 L_{f,g} \to \Omega^{\infty + 1} E(f,g)^\xi
$$
where $E(f,g)$ is the homotopy fiber product of the maps
$f$ and $g$ and $\xi$ is the virtual bundle or rank $n-p-q$
described as in the
introduction.

\begin{defn} The composite map
$$
{\mathfrak L}_{f,g}\: F_{f,g} \to L_{f,g} \to \Omega^{\infty + 1} E(f,g)^\xi
$$
is called the 
{\it generalized linking map} of the pair $(f,g)$.

The induced map on path components
$$
\pi_0({\mathfrak L}_{f,g})\: \pi_0(F_{f,g}) \to \Omega_{p+q+1-n}(E(f,g);\xi)
$$
are called the {\it generalized linking function} of $f$ and $g$.
\end{defn}

\begin{ex} Take $N = D^n$, $P = S^p$ and $Q = S^q$.
Take $f$ and $g$ to be the constant maps. 
Then $\pi_0(F_{f,g})$
is identified with the set of link homotopy classes of {\it link maps}
$$
F\: S^p \amalg S^q \to {\Bbb R}^n\, ,
$$ 
i.e., maps such that $F(S^p) \cap F(S^q) = \emptyset$.
This set is often denoted by
$$
L_{p,q}^n \, .
$$

The space $E(f,g)$ is identified with $S^p \times S^q$, 
and the virtual bundle $\xi$ is trivial. If we additionally
assume that $p,q \le  n-2$, then the
bordism group $\Omega_{p+q+1-n}(E(f,g);\xi)$ is isomorphic
to the stable homotopy group
$$
\pi_{p+q+1 -n}^{\text{\rm st}}(S^0)\, .
$$
With respect to these identifications, the 
generalized linking function takes the form
$$
\pi_0({\mathfrak L}_{f,g})\: L_{p,q}^n \to \pi_{p+q+1 -n}^{\text{\rm st}}(S^0)\, ,
$$
and is just Massey and Rolfsen's  ``$\alpha$-invariant''
(see \cite{Massey-Rolfsen} and \cite{Koschorke}).
When $p+q+1 = n$,  the invariant
is an integer and is just the classical Hopf linking number. 
\end{ex}

\begin{rem} When $p+q+1 =n$, 
the generalized linking function is closely related
to the affine linking invariants defined in a recent paper
of Chernov and Rudyak \cite{Chernov-Rudyak}. The
difference between our invariants and theirs has to
do with indeterminacy.

Namely, a point in $F_{f,g}$ consists of a pair
of maps $f_1\: P \to N$ and $g_1\: Q \to N$ having disjoint images, together
with a choice of homotopies from $f_1$ to $f$ and $g_1$ to $g$. The
choice of homotopies allows us to get an invariant living in the
full bordism group. 

By constrast, the Chernov and Rudyak do not choose the homotopies, 
but this comes at the cost of obtaining an invariant
landing in a certain {\it quotient} of the bordism group. The subgroup
which one quotients out by measures the indeterminacy arising 
from the choice the homotopies.
\end{rem}

\section{Fixed point theory}

Suppose we are given a closed manifold $M$ of dimension $m$
and a self map $f\:M \to M$. Then the {\it graph} of $f$
is the map
$$
G_f = (\text{id},f)\:M \to M \times M \, .
$$
Then $G_f$ is homotopic
to a map missing the diagonal $\Delta \subset
M \times M$, if and only if $f$ is homotopic to a fixed-point free map.

Actually, quite a bit more is true. 
Let $$\text{\rm map}^\flat(M,M) \subset \text{\rm map}(M,M)$$
be the subspace of fixed point free self maps of $M$.
Then we have a commutative square
\begin{equation}\label{fixed-point-sq} 
\SelectTips{cm}{}
\xymatrix{
\text{\rm map}^\flat(M,M) \ar[r] \ar[d] & \text{\rm map}(M,M) \ar[d]\\
\text{\rm map}(M,M^{\times 2}-\Delta) \ar[r] &\text{\rm map}(M,M^{\times 2})
}
\end{equation}
whose horizontal arrows are inclusions and whose vertical ones
are given by taking the graph of a self map. 
Since the first factor projection $M^{\times 2} -\Delta \to M$
is a fibration, the induced map of mapping spaces
$$
\text{\rm map}(M,M^{\times 2}-\Delta) \to 
\text{\rm map}(M,M)
$$
is also a fibration whose fiber at the identity is
$\text{\rm map}^\flat(M,M)$. Similarly, the
first factor projection map $M^{\times 2} \to M$
gives a fibration
$$
\text{\rm map}(M,M^{\times 2}) \to 
\text{\rm map}(M,M)
$$
whose fiber at the identity is $\text{\rm map}(M,M)$. Consequently, the
square \eqref{fixed-point-sq} is homotopy cartesian. 

We infer that the homotopy fiber of the map
$$
\text{\rm map}^\flat(M,M) \to \text{\rm map}(M,M)
$$
taken at $f$ coincides with the homotopy
fiber of the map 
$$
\text{\rm map}(M,M^{\times 2}-\Delta) \to \text{\rm map}(M,M^{\times 2})
$$ taken at $G_f$. 

Therefore, solutions of the intersection problem
$$
\SelectTips{cm}{}
\xymatrix{
& M^{\times 2} - \Delta \ar[d] \\
M \ar[r]_{G_f} \ar@{..>}[ur]
& M^{\times 2} \, 
}
$$
correspond to deformations of $f$ to a fixed point
free map. Furthermore, the relevant
virtual bundle $\xi$ is this case is trivial of
rank $0$.

Hence, applying
Theorems \ref{main} and \ref{families}, we immediately obtain

\begin{thm} \label{fixed-points}

{\flushleft (1).}  Let $f\: M \to M$ be a
self map of a closed manifold $M$.
Then is a well-defined class
$$
\ell(f) := \chi(G_f,\Delta) \in \Omega^{\text{\rm fr}}_0(L_f M)\, ,
$$
where $L_f M$ denotes the space of paths 
$\gamma\: [0,1]\to M$ such that $\gamma(1) = f(\gamma(0))$, and
$\Omega^{\text{\rm fr}}_0(L_f M)$ is its framed bordism group in degree
zero. When $f$ is homotopic to a fixed-point free map, then 
$\ell(f)$ is trivial.

Conversely, if $m \ge 3$, and $\ell(f)$ is trivial, then $f$ is homotopic to a
fixed-point free map.
\medskip

{\flushleft (2).} Assume $\ell(f)$ is trivial. 
Let ${\cal F}_f$ be denote the
homotopy fiber of the map 
$$
\text{\rm map}^\flat(M,M) \overset{\subset}\to \text{\rm map}(M,M)\, ,
$$
taken at $f$. Then there is an $(m{-}3)$-connected map
$$
{\cal F}_f \to \Omega Q(L_f M)_+ \, .
$$
where the target is the loop space of the stable homotopy
functor applied to $L_f M$ with the union of a disjoint basepoint.
\end{thm}

\begin{rems} (1). The map $L_f M \to \text{pt}$
induces a homomorphism
$$
\Omega^{\text{\rm fr}}_0(L_f M) \to \Omega^{\text{\rm fr}}_0(\text{pt}) 
\cong {\Bbb Z}\, .
$$
It is possible to show that image of
$\ell(f)$ under this homomorphism is just the 
classical {\it Lefschetz trace} of $f$.

In fact, the group $\Omega^{\text{\rm fr}}_0(L_f M)$ 
is isomorphic to the coinvariants 
of $\pi = \pi_1(M,*)$ acting on its integral group ring
by means of {\it twisted conjugation:} $(g,x) \mapsto 
g x \phi (g)^{-1}$, for $g,x \in \pi$
where $\phi\: \pi \to \pi$ is the homomorphism
induced by $f$ and  a choice of path from $*$ to $f(*)$.

Using our Index Theorem (\ref{index}),
 $\ell(f)$ can be identified with the 
{\it Reidemeister trace of $f$} \cite{Reidemeister}. 
With respect to  this identification, the first part of \ref{fixed-points}
is essentially equivalent to the homotopy
converse to the Lefschetz fixed point theorem.
(cf.\ \cite{Wecken}, \cite{Brown}, \cite[Cor.\ 3]{Brown_book})
\smallskip

{\flushleft (2).} The self map $f$ describes an action of the natural
numbers $\Bbb N$ on $M$. The space $L_f M$ is just the
homotopy fixed point set $M^{h\Bbb N}$ of this action, and
the bordism class of the theorem is the one associated with the
evident inclusion
$$
M^{\Bbb N} \subset M^{h\Bbb N}
$$
from fixed points to homotopy fixed points.
\smallskip

{\flushleft (3).} The space $Q(L_f M)_+$ has an alternative description as
the  {\it topological
Hochschild homology} of the associative ring (spectrum) $$S^0[\Omega M]$$ 
(= the suspension spectrum of $(\Omega M)_+$)
with coefficients in the bimodule given by $S^0[\Omega M]$,
where the action is defined by twisted conjugation of loops.
\end{rems}

\subsection*{Families over a disk}
$$
F\: M\times D^j \to M
$$
be a $j$-parameter family of self maps of $M$ such that
the restriction of $F$ to $M \times S^{j-1}$ is fixed point free.
Let $f\: M \to M$ be the the self map associated with the basepoint
of $S^{j-1}$.

Then $F$ defines a map of pairs 
$$
(D^j,S^{j-1}) \to (\maps(M,M),\maps^\flat(M,M))
$$
which by \ref{fixed-points} yields an invariant 
$$
\ell^{\text{fam}}(F) \in \Omega_j^{\rm fr}(L_f M)\, .
$$

\begin{cor} Assume $j \ge 1$.
If $F\:M \times D^j \to M$ is homotopic
rel $M \times S^{j-1}$ to a fixed point free family,
the $\ell^{\text{\rm fam}}(F) = 0$.

If $j \le m-3$ and  $\ell^{\text{\rm fam}}(F) = 0$,
then the converse also holds.
\end{cor}

\begin{rem} Using different methods,
an invariant of this type was discovered
by  Geoghegan and Nicas (see \cite{G-N}), who
have intensively studied the $j = 1$ case.
\end{rem}

\subsection*{General families} 
The previous situation can be generalized as follows: suppose
that $p\: E \to B$ is a smooth fiber bundle over a connected space $B$
having closed manifold fibers.
Given a fiberwise self map $f\: E \to E$, one asks whether
$f$ is fiberwise homotopic to a map having no fixed points
(the case when $B$ is disk coincides with the previous
situation). This type of fixed point problem  was studied by 
Dold \cite{Dold}.

For $b\in B$, let $f_b \: E_b \to  E_b$
be the induced self map of the fiber at $b$. Then we
have a fibration 
$$
q\:{\cal E}_f \to B
$$
whose fiber at $b$ is $Q (L_{f_b}E_b)_+$.
Note that the space of sections of $q$ is an
infinite loop space. So if we 
let
$$
\Omega^{B,\text{\rm fr}}_0({\cal E}_f)
$$
denote the set of homotopy classes of sections of $q$, 
we see that the latter has the structure  of an abelian group.

\begin{thm} \label{twisted} There is a well-defined element
$$
{\ell}_B(f) \in \Omega^{B,\text{\rm fr}}_0({\cal E}_f)\, ,
$$
which vanishes when $f$ is fiberwise homotopic to a
fixed point free self map.

Conversely, assume the fibers of $p$ have dimension $m$,
and $B$ is homotopically the retract of a cell complex with
cells of dimension $\le m - 3$. Then ${\ell}_B(f) = 0$ implies
that $f$ is fiberwise homotopic to a fixed point free self map.
\end{thm}

A version of this theorem was proved by V.\ Coufal \cite{Coufal},
relying on the fibered theory of 
Hatcher and Quinn \cite[th.\ 4.2]{H-Q}.

We will not give the 
proof of \ref{twisted} here because a substantial revision
of our methods would be needed. Rather, we will be content
to explain the main issue.
The situation at hand is represented by a family of
intersection problems parametrized by points of some
parameter space $B$.
Hence, \ref{twisted} would follow from a suitable generalization
of Theorem \ref{families}. 

The required generalization is this: we are given
a parameter space $B$ and a (continuous) family of intersection problems
$$
\SelectTips{cm}{}
\xymatrix{
& N_b - Q_b \ar[d] \\
P_b \ar[r]_{f_b}\ar@{..>}[ur]
& N_b \, 
}
$$
parametrized by $b \in B$. The approach of this paper is
to {\it choose a basepoint} in each $P_b$. 
Once a continuously varying family
of basepoints is chosen, the invariant can be defined,
and furthermore, the parametrized intersection problem 
can be attacked using the machinery which proves Theorem \ref{families}.
The problem with this approach is: {\it a 
continuously varying family of basepoints might not exist.} 

The reason we needed to choose basepoints is that we converted
a {\it fiberwise} problem---finding a section of some fibration---into an {\it equivariant} one, in which the group
is a loop space of the base space of the fibration
(the conversion was made by taking a thick fiber).
The basepoint was needed to form the loop space.

Hence, to avoid basepoint issues, we should not have
passed to the equivariant category at all. That is, we should
have developed all the technology in the fiberwise setting.
A review of our constructions shows that the crux of the matter
is to reformulate Theorem \ref{duality}. Recall
that this result relies heavily on the 
main results of \cite{Klein_dualizing}, and the latter paper
is written in the equivariant context.
Fortunately, the results of that paper 
can indeed be redone in the fiberwise setting (see e.g., \cite{PoHu}).  
Rather than deluge the reader with
technical details, we will leave it
to him/her to transform this discussion into
a proof of \ref{twisted}.

\begin{rem} Considering the map $E \to B$ as a map
of spaces over $B$, we get an induced homomorphism
$$
\Omega^{B,\text{fr}}_0 ({\cal E}_f) \to \Omega^{B,\text{fr}}_0 (\text{pt}) \cong\pi^0_{\text{st}}(B_+)
$$
where the target is the stable cohomotopy of $B$ in degree zero.

Let $\text{map}_B(E,E)$ denote the space of fibered self maps of 
$E$. Then \ref{twisted}, followed by above, yields a function 
$$
I\:\pi_0(\text{map}_B(E,E)) \to \pi^0_{\text{st}}(B_+)\, ,
$$
which presumably coincides with Dold's fixed point index
\cite{Dold} (compare \cite{Crabb-James}). 
\end{rem}

\section{Disjunction}

A {\it disjunction problem} is an intersection problem
$$
\SelectTips{cm}{}
\xymatrix{
& N - Q \ar[d] \\
P \ar[r]_{f}\ar@{..>}[ur]
& N \, 
}
$$
in which $f\: P \to N$ is a smooth embedding.
In this instance we ask for an isotopy of $f$ to another
embedding $g$ such that the image of $g$ is disjoint from $Q$.
Let $$\text{emb}(P,N)$$ be the space of embeddings of $P$ in $N$.

The following is a special instance
of the main results of \cite{G-K}. Its proof
involves a mixture of homotopy theory, surgery, and Morlet's
disjunction lemma (\cite{Morlet}, \cite{BLR}):

\begin{thm}[Goodwillie-Klein \cite{G-K}] \label{disjunction}
Assume $p,q \le n-3$. Then the commutative square
$$
\SelectTips{cm}{}
\xymatrix{
\text{\rm emb}(P,N-Q) \ar[r] \ar[d] & \text{\rm emb}(P,N)\ar[d] \\
\maps(P,N-Q) \ar[r]  & \maps(P,N)
}
$$
is $(n-2p-q-3)$-cartesian.
\end{thm}

{\flushleft (A square is {\it $j$-cartesian} if the map from its initial
term to the homotopy pullback of the other terms is $j$-connected.)}

\begin{rem} Hatcher and Quinn state a version of \ref{disjunction}
(cf.\ \cite[th.\ 4.1]{H-Q}), but give a proof which is
only a few lines long, and is therefore too scant 
to be regarded as complete. Additionally, we expected to find a parametrized
``concordance implies isotopy'' problem in their argument, but nothing
of the sort is mentioned. 
\end{rem}

\begin{cor}[Hatcher-Quinn {\cite[th.\ 1.1]{H-Q}}] Assume $2p+q+3\le n$.
Then an embedding $f\: P \to N$
has a isotopy to an embedding with image disjoint from $Q$
if and only if there is a homotopy of $f$ to a map
having image disjoint from $Q$.
\end{cor}

\begin{cor}[compare {\cite[th.\ 2.2]{H-Q}}] Assume $p+2q,2p + q \le n-3$.
Then $f$ is isotopic to an embedding with image disjoint
from $Q$ if and only if $\chi(f,i_Q) = 0$.
\end{cor}

Let ${\cal E}_{f}$ be the homotopy fiber at $f$ of 
$$
\text{emb}(P,N-Q) \to \text{emb}(P,N)\, .
$$
Then \ref{disjunction} gives an $(n-2p-q-3)$-connected map
$$
{\cal E}_{f} \to {\cal F}_f\, .
$$
Compose this with the map of Theorem \ref{families}
to get a map 
$$
{\cal E}_{f}\to \Omega^{\infty +1} E(f,i_Q)^\xi\, .
$$

\begin{cor} The map ${\cal E}_{f}\to \Omega^{\infty +1} E(f,i_Q)^\xi$
is $$\min(2n-2q-p-3,2n-2p-q-3){\text{\rm -connected.}}$$
\end{cor}

\section{An index theorem}

\subsection*{The Hatcher-Quinn invariant}
Given an intersection problem
$$
\SelectTips{cm}{}
\xymatrix{
& N - Q \ar[d] \\
P \ar[r]_f\ar@{..>}[ur]
& N \, ,
}
$$
we will assume $f$ is transverse to $i_Q \: Q \subset N$. Then
the intersection manifold
$$
D := f^{-1}(Q) 
$$
has dimension $p+q-n$ and comes equipped with
a map $D \to E(f,i_Q)$ given by $x \mapsto (x,c_{f(x)},f(x))$, where $c_{f(x)}$
is the constant path. The pullback of $\xi$ along this map coincides
with the stable normal bundle of $D$. The bordism class of these
data in $\Omega_{p+q-n}(E(f,i_Q);\xi)$ is called the {\it Hatcher-Quinn
invariant}, and is denoted
$$
[f \pitchfork i_Q]\, .
$$
Recall that $E\to P$ is the fibration given by 
taking the homotopy pullback of the inclusion $N-Q \to N$ along
the map $f\: P \to N$. It has a stable cohomotopy Euler class
$e(p) \in \pi_0(\Sec^{\rm st}(S_P E \to P))$.

\begin{thm}[``Index Theorem''] \label{index} There is a bijection
$$
{\mathfrak p\mathfrak d}\: \pi_0(\Sec^{\rm st}(S_P E \to P)) 
\overset \cong \to
\Omega_{p+q-n}(E(f,i_Q);\xi)
$$
such that $$\mathfrak p \mathfrak d (e(p)) \,\, = \,\, [f\pitchfork i_Q]\, .$$
\end{thm}

\begin{proof}
Consider the fiberwise suspended fibration $S_P E \to P$. It is identified
with the homotopy pullback of 
$$
S_N (N-Q) \to N
$$
along $f\: P \to N$. Consequently, specifying a section of $S_P E \to P$
is equivalent to choosing a map $s\:P \to S_N(N-Q)$ and also a homotopy from
the composite
$$
P \overset s \to S_N (N-Q) \to N
$$
to $f$.

Let $\nu$ be the normal bundle of $Q$ in $N$,
and identify its unit disk bundle $D(\nu)$ with a tubular neighborhood
of $Q$. Let $C$ denote the closure of the complement of this neighborhood.
The we have a homeomorphism
$$
N \cup_{S(\nu)} D(\nu) \cong N \cup_C N \, .
$$
Furthermore, $N \cup_C N$ is identified with
$S_N (N-Q)$ up to weak equivalence. 

Let $\pi\: N\cup_{S(\nu)} D(\nu) \to N$ be the evident map.
Then up to homotopy, choosing a section of $S_P E \to P$ is
the same thing as specifying a map
$$
\sigma\: P \to N \cup_{S(\nu)} D(\nu)
$$
and a homotopy from $\pi\sigma$ to $f$. 

Let
$z_0\: Q \to D(\nu)$ denote the zero section.
If $\sigma$ is transverse to $z_0$ then 
$$
{\cal D} := P \cap \sigma^{-1}(z_0(Q))
$$
is a submanifold of $P$ of dimension $p+q-n$ which comes
equipped with a map  ${\cal D} \to E(f,i_Q)$. The normal bundle
of ${\cal D}$ in $P$ is isomorphic to the pullback of $\nu$. Therefore
the stable normal bundle of ${\cal D}$ 
is a Whitney sum of the pullbacks of $\nu$ and
the stable normal bundle of $P$. 
This implies that the stable normal bundle of ${\cal D}$ is given by 
taking the pullback of $\xi$. Consequently, these data determine an element of 
$\Omega_{p+q-n}(E(f,i_Q);\xi)$. 

If $\sigma$ isn't transverse to $z_0$, we can still 
make it transverse without
altering its vertical homotopy class. 
Therefore, we get  a well-defined function
$$
\pi_0(\Sec(S_P E \to P)) \to \Omega_{p+q-n}(E(f,i_Q);\xi)\, .
$$
By definition, the vertical homotopy class of 
the section $s_+$ maps to
$[f\pitchfork i_Q]$.

If we are now given a stable section of 
$S_P E \to P$, the above construction is now applied
to an associated map 
$$
\sigma \: P \times D^j \to N \cup_{S(\nu\oplus\epsilon^j)} 
D(\nu\oplus \epsilon^j) 
$$
for sufficiently large $j$. This produces 
a submanifold ${\cal D} \subset P \times D^j$ together with 
a map ${\cal D}\to E(f,i_Q)$ compatible with the bundle data. So we get
a function
$$
{\mathfrak p \mathfrak d}\:
\pi_0(\Secst(S_P E \to P)) \to \Omega_{p+q-n}(E(f,i_Q);\xi)\, .
$$
which maps $[s_+]$ to $[f\pitchfork i_Q]$.

We now show
that  $\mathfrak p \mathfrak d$ is onto. Let $(V,g,\phi)$ represent
an element of  $\Omega_{p+q-n}(E(f,i_Q);\xi)$. Assume first that
the composite
$$
V \overset g \to E(f,i_Q) \to P
$$
is an embedding. Notice that $\phi$ provides a stable isomorphism
of the normal bundle $\nu_{V \subset P}$ of this embedding with the 
pullback of normal bundle $\nu$ of $Q \subset N$ along $f$.
For the moment, we also 
assume that $\phi$ arises from an unstable isomorphism of
these bundles. By choosing a tubular
neighborhood, we get an embedding $D(\nu_{V \subset P}) \subset P$.
Let $C$ denote the closure of the complement of this neigbhorhood. Then
we have a decomposition
$$
P \,\, = \,\, C\cup D(\nu_{V \subset P}) \, ,
$$
so we get a map 
$$
\begin{CD}
 C\cup D(\nu_{V \subset P}) @> f \cup \phi >> N \cup D(\nu)
\end{CD} \, .
$$
By the remarks at the beginning of the proof, we see that this
determines a homotopy class of section of $S_P E \to P$.

In the general case, we choose an embedding of $V$ in $D^j$. Then
$V$ is embedded in the product $P \times D^j$, and if $j$ is sufficiently
large the normal bundle of this embedding is identified with the pullback
of the $\nu\oplus \epsilon^j$ using $\phi$. 
The construction above now gives a homotopy class
of stable section of $S_P E \to P$. It is straightforward to check that 
$\mathfrak p \mathfrak d$ applied to this homotopy class recovers the original
bordism class. Consequently, $\mathfrak p \mathfrak d$ is onto.

The proof that $\mathfrak p \mathfrak d$ is one-to-one is similar. In this
case however, one starts with a bordism whose boundary arises from two 
stable sections of $S_P E \to P$. Then, working relatively, one constructs
a homotopy between these stable sections. We leave these details to the reader.
This completes the proof of Theorem \ref{index}.
\end{proof}

\end{document}